\documentclass[12pt,reqno]{amsart}
\usepackage{amsmath}
\usepackage{amssymb}
\usepackage{amsthm}
\usepackage{mathrsfs}
\usepackage{latexsym}
\usepackage{xcolor}
\usepackage{comment}

\newtheorem{theorem}{Theorem}[section]
\newtheorem{lemma}{Lemma}[section]
\newtheorem{corollary}{Corollary}[section]
\newtheorem{remark}{Remark}[section]
\newtheorem{definition}{Definition}[section]
\newtheorem{proposition}{Proposition}[section]
\newtheorem{example}{Example}[section]
\newtheorem{assumption}{Assumption}[section]
\numberwithin{equation}{section}
\newcommand{\bth}{\begin{theorem}}
\newcommand{\ethe}{\end{theorem}}

\newcommand{\bre}{\begin{remark}}
\newcommand{\ere}{\end{remark}}

\newcommand{\ble}{\begin{lemma}}
\newcommand{\ele}{\end{lemma}}
\newcommand{\bde}{\begin{definition}}
\newcommand{\ede}{\end{definition}}

\newcommand{\bco}{\begin{corollary}}
\newcommand{\eco}{\end{corollary}}

\newcommand{\bpr}{\begin{proposition}}
\newcommand{\epr}{\end{proposition}}

\newcommand{\bexer}{\begin{exercise}}
\newcommand{\eexer}{\end{exercise}}
\newcommand{\breh}{\begin{hint}}
\newcommand{\ereh}{\end{hint}}

\newcommand{\halmos}{\hfill \qed}

\newcommand{\bexam}{\begin{example}}
\newcommand{\eexam}{\end{example}}

\newcommand{\pr} {{\bf Proof.}}

\newcommand{\bfi}{\begin{fig}}
\newcommand{\efi}{\end{fig}}

\newcommand{\beao}{\begin{eqnarray*}}
\newcommand{\eeao}{\end{eqnarray*}\noindent}

\newcommand{\beam}{\begin{eqnarray}}
\newcommand{\eeam}{\end{eqnarray}\noindent}
\newcommand{\E}{\mathbf{E}}
\newcommand{\PP}{\mathbf{P}}

\newcommand{\xto}{x\to\infty}

\newcommand{\bF}{\overline{F}}

\newcommand{\bV}{\overline{V}}

\newcommand{\bbr}{{\mathbb R}}

\newcommand{\bbb}{{\mathbb B}}

\newcommand{\bbn}{{\mathbb N}}

\newcommand{\vep}{\varepsilon}

\allowdisplaybreaks[1]
\begin{document}
\title[Uniform asymptotics with multivariate subexponential claims]{Uniform asymptotics for a multidimensional renewal risk model with multivariate subexponential claims}

\author[D.G. Konstantinides, J. Liu, C.D. Passalidis]{Dimitrios G. Konstantinides,\,JiaJun Liu, Charalampos  D. Passalidis} 

\address{Dept. of Statistics and Actuarial-Financial Mathematics,
University of the Aegean,
Karlovassi, GR-83 200 Samos, Greece}
\email{konstant@aegean.gr}
\address{Dept. of Financial and Actuarial Mathematics,
Xi'an Jiaotong-Liverpool University,
Suzhou, Jiangsu, China}
\email{JIAJUN.LIU@xjtlu.edu.cn}
\address{Dept. of Statistics and Actuarial-Financial Mathematics,
University of the Aegean,
Karlovassi, GR-83 200 Samos, Greece}
\email{sasd24009@sas.aegean.gr.}

\date{{\small \today}}

\begin{abstract}
In this paper, we investigate a multidimensional risk model driven by a common renewal process under a constant force of interest. The claim sizes generated by each line of business are independent and identically distributed random vectors with possibly dependent components, and their common distribution belongs to the class of multivariate subexponential distributions. We derive locally uniform asymptotic estimates for the probability that discounted aggregate claims enter certain `rare sets', and extend these to uniform estimates over all time horizons under some extra mild conditions. As a direct application, we obtain uniform estimates for the finite-time ruin probability defined using various ruin sets. Additionally,  we provide examples of distributions belonging to these multivariate heavy-tailed classes, which are not limited to the case of multivariate regular variation.
\end{abstract}

\maketitle
\textit{Keywords: multidimensional risk model; renewal process; multivariate subexponentiality; multivariate positively decreasing distributions; uniformity}
\vspace{3mm}

\textit{Mathematics Subject Classification}: Primary 62P05 ;\quad Secondary 60G70.


\section{Introduction} \label{sec.KLP.1}

In this paper, we consider an insurer operating $d$-lines of business $(d \in \bbn)$
that share a common counting process for claim vectors $\{{\bf X}^{(i)},\,i \in \bbn\}.$ These claim vectors form a sequence of independent and identically distributed (i.i.d.) non-negative random vectors with a common distribution $ F.$ 

Note that each claim vector $\mathbf{X}^{(i)} = \left(X_1^{(i)}, \ldots, X_d^{(i)}\right)$, for $i \in \mathbb{N}$, may contain zero components (but not all zero), and arrives at time $\tau_i$, with $\tau_0 = 0$ by convention. The sequence $\{\tau_i\}_{i \in \mathbb{N}}$ forms a renewal counting process
\beao
N(t) :=\sup\{n\in \bbn\;:\;\tau_n \leq t\}\,,
\eeao
namely, the inter-arrival times $\{\theta_i\}_{i \in \mathbb{N}}$, defined by $\theta_i = \tau_i - \tau_{i-1}$ for $i \geq 1$, constitute a sequence of i.i.d.\ non-negative and non-degenerate to zero random variables. The renewal process $\{N(t), \, t \geq 0\}$ is assumed to have a {finite renewal function}:
\beao
\lambda(t) = \E[N(t)] = \sum_{i=1}^{\infty} \PP[\tau_i \leq t]\,.
\eeao

For applicability to risk theory, we restrict our analysis to the interval
$\Lambda:=\{t\;:
\lambda(t)>0 \}= \{t\;:\;\PP[\tau_1 \leq t] >0\}$.
This restriction is natural in both univariate and multivariate risk models, as the ruin probability trivially equals zero if there is no positive probability of claim occurrence.

Assuming the insurer invests surplus at a constant interest force $r \geq 0$ (with $r = 0$ representing no investment), the {discounted aggregate claims} at time $t \geq 0$ are given by
\beam \label{eq.KLP.1.2}
{\bf D}_r(t) =\sum_{i=1}^{N(t)} {\bf X}^{(i)}\,e^{-r\,\tau_i}=\left( 
\begin{array}{c}
\sum_{i=1}^{N(t)} X_{1}^{(i)}\,e^{-r\,\tau_{i} } \\ 
\vdots \\ 
\sum_{i=1}^{N(t)} X_{d}^{(i)}\,e^{-r\,\tau_{i} } 
\end{array} 
\right)\,.
\eeam

Furthermore, we consider that the initial capital of the insurer is equal to $x>0$ and is allocated through the deterministic weights $l_1,\,\ldots,\,l_d > 0$, with 
\beao
\sum_{i=1}^d l_i =1\,,
\eeao
into $d$-lines of business. The insurer receives premiums at a rate given by the vector ${\bf p}(t)=(p_1(t),\,\ldots,\, p_d(t))$, for any $t \geq 0$, where each premium density $p_i(t)$ is non-negative and bounded above by a positive constant  $C_i\in(0,\infty)$, for all $i=1,\,\ldots,\,d$ and $t\geq 0$. This technical assumption of bounded premium densities is a mild condition, since the $C_i$ can be chosen arbitrarily large. This flexibility also permits arbitrary dependence between the premium process and other sources of randomness. Therefore, the discounted surplus process at time $t \geq 0$ is given by the relation
\beam \label{eq.KLP.1.3}
{\bf U}(t):=\left( 
\begin{array}{c}
U_{1}(t) \\ 
\vdots \\ 
U_{d}(t) 
\end{array} 
\right) =x\,\left( 
\begin{array}{c}
l_{1} \\ 
\vdots \\ 
l_{d} 
\end{array} 
\right) +\left( 
\begin{array}{c}
\int_{0-}^{t}e^{-r\,s}\,p_1(s)\,ds \\ 
\vdots \\ 
\int_{0-}^{t}e^{-r\,s}\,p_d(s)\,ds 
\end{array} 
\right) -{\bf D}_r(t)\,.
\eeam

In modern insurance, the need to maintain competitive multi-risky portfolios has generated considerable research interest in the multidimensional risk model \eqref{eq.KLP.1.3}. For the special case $d=2,$ numerous studies have analyzed this framework under varying distributional assumptions and dependence structures for claims, including \cite{chen:wang:wang:2013}, \cite{gao:yang:2014}, \cite{jiang:wang:chen:xu:2015}, \cite{yang:wang:liu:zhang:2019}, \cite{chen:li:cheng:2023}, \cite{chen:cheng:zheng:2025}, and \cite{xu:shen:wang:2025}, among others.

For higher-dimensional settings, several studies have investigated the risk model in \eqref{eq.KLP.1.3} under the assumption that claim vectors exhibit multivariate regular variation (MRV). Key examples include \cite{hult:lindskog:2006}, \cite{konstantinides:li:2016}, \cite{li:2016}, \cite{cheng:konstantinides:wang:2022}, \cite{yang:su:2023}, and \cite{cheng:konstantinides:wang:2024}. Most of these studies additionally presume that the components within each claim vector are asymptotically dependent, typically enforced through an MRV structure. This assumption necessitates that the marginal distributions fall within a regularly varying class, thereby constraining the framework to scenarios characterized by sufficiently extreme events. Consequently, numerous moderately heavy-tailed distributions, such as the lognormal and Weibull distributions, are excluded from consideration.

Based on this gap in the literature, following the pioneering paper \cite{samorodnitsky:sun:2016}, we study the model in \eqref{eq.KLP.1.3} under some multivariate distribution classes larger than MRV, and under arbitrarily dependent components of each claim vector. See related models in \cite{samorodnitsky:sun:2016},  \cite{konstantinides:passalidis:2024h}, \cite{konstantinides:passalidis:2024j}, \cite{passalidis:2025k}. In this paper, we follow the claim-vectors model, with distribution classes larger than MRV and with arbitrarily dependent components. We also establish uniform asymptotic estimations for the behavior of the entrance probability of the discounted aggregate claims into some `rare sets', and further for the ruin probability of the insurer.

Namely, we establish the asymptotic relation
\beam \label{eq.KLP.1.4}
\PP\left[ {\bf D}_r(t) \in x\,A  \right] \sim \int_0^t \PP[{\bf X} \in x\,e^{r\,s}\,A] \,\lambda(ds)\,,
\eeam
as $\xto$, uniformly for any $t \in \Lambda_T :=\Lambda\cap (0,\,T]$, for some fixed constant $T \in \Lambda$, when the distribution of the claim-vectors belongs to the class of multivariate subexponential distributions, see Section 2. The rare set $x\,A$ can represent various actuarially relevant scenarios. In the univariate case where $A = (1,\infty)$, relation \eqref{eq.KLP.1.4} reduces to the classical asymptotic framework for ruin probabilities in one-dimensional risk models, as seen in \cite{tang:2007}, \cite{hao:tang:2008}, \cite{li:2017}, among others.

Under a slightly narrower distributional assumption, specifically, the class of multivariate positively-decreasing subexponential distributions (see Section 2 for the definition), we strengthen the local uniformity condition in \eqref{eq.KLP.1.4} to a global uniformity condition, ensuring that \eqref{eq.KLP.1.4} holds uniformly for all $t \in \Lambda.$ A direct consequence of these results is a uniform asymptotic estimate for the finite-time ruin probability, where ruin is defined as the surplus entering a predefined ‘ruin set’.

The rest of the paper is organized as follows. In Section 2, we provide the preliminaries on multivariate distribution classes with heavy tails. In Section 3, we present our main results. In Section 4, we construct some distribution examples that belong to the classes of multivariate subexponential and multivariate positively decreasing subexponential distributions, primarily based on a concrete `rare set'. In Sections 5 and 6, we present the proofs of the main results, along with some necessary lemmas.

\section{Preliminaries} \label{sec.KLP.2}

In this section, we present some preliminary concepts for multivariate heavy-tailed distributions that will be used later.

\subsection{Notation} \label{sec.KLP.2.0}

In what follows, all the limit relations hold as $\xto$ unless otherwise stated. All the vectors are denoted with bold script, and the dimension is equal to $d$, for some $d\in \bbn$. Namely, ${\bf x} =(x_1,\,\ldots,\,x_d)$, and with ${\bf x}^T$ we denote the transpose of vector ${\bf x},$ and conventionally, the ${\bf 0} =(0,\,\ldots,\,0)$ represents the origin of the axes. For any two vectors ${\bf x},\,{\bf y} $, the component-wise sum and subtraction are denoted by ${\bf x}\pm{\bf y} =(x_1\pm y_1,\,\ldots,\,x_d\pm y_d),$ and the scalar product, with some positive quantity $k>0$, is denoted by $k\,{\bf x} =(k\,x_1,\,\ldots,\,k\,x_d)$. Furthermore, for a set $\bbb$ from space $\bbr^d$, we denote by $\bbb^c$ its complement set, by $\partial \bbb$ its boundary, and by $\overline{\bbb}$ its closed hull. For two $d$-dimensional positive functions ${\bf f},\,{\bf g}$ and some set $\bbb \subset \bbr_+^d \setminus \{{\bf 0}\}=[0,\,\infty)^d \setminus \{{\bf 0}\}$, we denote by ${\bf f}(x\,\bbb) =o[ {\bf g}(x\,\bbb)]$ if it holds 
\beao
\lim \dfrac{{\bf f}(x\,\bbb)}{{\bf g}(x\,\bbb)}=0\,,
\eeao
and we denote by ${\bf f}(x\,\bbb) \lesssim {\bf g}(x\,\bbb)$ (or equivalently ${\bf g}(x\,\bbb) \gtrsim {\bf f}(x\,\bbb) $) if it holds
\beao
\limsup \dfrac{{\bf f}(x\,\bbb)}{{\bf g}(x\,\bbb)}\leq 1\,,
\eeao
and we denote by ${\bf f}(x\,\bbb) \asymp {\bf g}(x\,\bbb)$, if it holds
\beao
0<\liminf \dfrac{{\bf f}(x\,\bbb)}{{\bf g}(x\,\bbb)}\leq \limsup \dfrac{{\bf f}(x\,\bbb)}{{\bf g}(x\,\bbb)}< \infty\,.
\eeao
For the $(d+1)$-variate, positive functions ${\bf f}^*$, ${\bf g}^*$, we say that it holds ${\bf f}^*(x\,\bbb,\,t) \sim {\bf g}^*(x\,\bbb,\,t)$, uniformly for any $t\in E$, with some set $E \neq \emptyset$, if
\beao
\lim_{\xto} \sup_{t\in E} \left| \dfrac{{\bf f}^*(x\,\bbb,\,t)}{{\bf g}^*(x\,\bbb,\,t)} -1\right| =0\,,
\eeao
and it holds  ${\bf f}^*(x\,\bbb,\,t) \lesssim {\bf g}^*(x\,\bbb,\,t)$ uniformly for any $t\in E$, if
\beao
\limsup_{\xto} \sup_{t\in E} \dfrac{{\bf f}^*(x\,\bbb,\,t)}{{\bf g}^*(x\,\bbb,\,t)} \leq 1\,.
\eeao

Finally, for a distribution $V$, we denote its tail by $\bV(x) = 1- V(x)$, for any $x\in \bbr$. In what follows, the random vectors have support in the non-negative quadrant $\bbr_+^d$. In this section, all the univariate distributions $V$ are such that $\bV(x)>0$, for any $x \in \bbr$.

\subsection{Multivariate distribution classes} \label{sec.KLP.2.1}

To present the definition of multivariate distribution classes, we need the following family of sets
\beao
\mathscr{R}:= \{A \subsetneq \bbr^d\;:\;A\, {\text open,\;increasing},\; A^c\; {\text convex}, \;{\bf 0} \notin \overline{A} \}\,,
\eeao
where the set $A$ is called increasing if, for any ${\bf x} \in A$ and ${\bf y} \in \bbr_+^d$, it holds ${\bf x}+{\bf y} \in A$. Let $\mathbf{X}$ be a random vector and $A \subsetneq \mathscr{R}$ be a fixed set. 
The random variable $Y_A$ is defined by
\beao
Y_A :=\sup \{ u\;:\; {\bf X} \in u\,A\}\,,
\eeao
and it has a proper distribution, denoted by $F_A$. 
The tail of $F_A$ is given by
\beao
\bF_A(x) :=\PP[{\bf X} \in x\,A] = \PP\left[ \sup_{{\bf p} \in I_A} {\bf p} ^T\,{\bf X} >x \right]\,,
\eeao
where $I_A \subsetneq \mathbb{R}^d$ is an index set for $A$, whose existence is established in \cite[Lem. 4.3(c)]{samorodnitsky:sun:2016}. 
The fact that $F_A$ is a proper distribution is verified in \cite[Lem. 4.5]{samorodnitsky:sun:2016}.

Hence, for any fixed set $A \subsetneq \mathscr{R}$, by \cite{samorodnitsky:sun:2016} we define the multivariate subexponentiality as follows. We say that the distribution $F$ belongs to the class $\mathcal{S}_A$, which means the class of multivariate subexponential distributions on the set $A$, symbolically $F \in \mathcal{S}_A$, if  $F_A \in \mathcal{S}$, namely if for any (or, equivalently, for some) integer $n\geq 2$, it holds
\beao
\lim \dfrac{\overline{F_A^{n*}}(x)}{\bF_A(x)} =n\,,
\eeao
whereby $F_A^{n*}$ denotes the $n$-th order convolution power of $F_A$. Additionally, multivariate subexponentiality over the whole set $\mathscr{R}$ was defined as class $\mathcal{S}_\mathscr{R} := \bigcap_{A\subsetneq \mathscr{R} } \mathcal{S}_A$.

On similar lines in \cite{konstantinides:passalidis:2024g}, we can find the class of multivariate long-tailed distributions. We say that the distribution $F$ follows a multivariate long-tailed distribution on $A$, symbolically $F \in \mathcal{L}_A$, if $F_A \in \mathcal{L}$, namely if for any (or, equivalently, for some) $a>0$ it holds
\beao
\lim \dfrac{\bF_A(x-a)}{\bF_A(x)} =1\,.
\eeao

For further details on the univariate classes $\mathcal{S}$ and $\mathcal{L}$, we refer to \cite{embrechts:klueppelberg:mikosch:1997} and \cite{leipus:siaulys:konstantinides:2023}. The class of multivariate positively decreasing distributions on a set $A$, denoted by $F \in \mathcal{P_D}_A$, was introduced in \cite{konstantinides:passalidis:2024h}. This class is characterized by distributions satisfying $F_A \in \mathcal{P_D}$, that is, for any (or, equivalently, for some) $v>1,$
\beao
\limsup \dfrac{\bF_A(v\,x)}{\bF_A(x)} <1\,,
\eeao
Similarly, the class of multivariate subexponential, positively decreasing distributions on $A$, denoted by $F \in \mathcal{A}_A$, is defined by $F_A \in \mathcal{A}$, where $\mathcal{A}:=\mathcal{S}\cap \mathcal{P_D}$. For the classes mentioned above, we adopt the general notation  $\mathcal{B}_\mathscr{R} :=\bigcap_{A \subsetneq \mathscr{R} } \mathcal{B}_A$, where $\mathcal{B} \in \{\mathcal{L},\,\mathcal{P_D},\,\mathcal{A},\,\mathcal{S} \}$. For the classes mentioned above, we adopt the general notation $\mathcal{P_D}$ is sufficiently broad to include distributions with both light and heavy tails. In fact, among all distributions with an infinite right endpoint, only the extended slowly varying distributions do not belong to $\mathcal{P_D}$ (\cite[sec. 7.4]{konstantinides:2018}). Consequently, most subexponential distributions belong to $\mathcal{A}$, which implies that the class $\mathcal{S}\setminus \mathcal{P_D}$ is relatively sparse.

The restriction to class $\mathcal{A}$ resolves several problems where class $\mathcal{S}$ fails, particularly concerning infinite sums, infinite-time ruin probabilities, and globally uniform asymptotic estimates (see, e.g., \cite{konstantinides:tang:tsitsiashvili:2002}, \cite{chen:ng:tang:2005}, \cite{hao:tang:2008} among others). The class $\mathcal{P_D}$ is characterized by the lower Matuszewska index $J_{F_A}^{-}$. Concretely, $F_A \in \mathcal{P_D}$ if and only if $J_{F_A}^{-}>0$, where
\beao
J_{F_A}^{-} =-\lim_{v\to \infty} \dfrac{\log \bF_A^*(v)}{\log v}\,,
\eeao
with 
\beao
\bF_A^*(v)= \limsup \dfrac{\bF_A(v\,x)}{\bF_A(x)}\,.
\eeao
For further results on Matuszewska indices, please see, e.g., \cite{bingham:goldie:teugels:1987}, \cite{bardoutsos:konstantinides:2011}, \cite{leipus:siaulys:konstantinides:2023}, and \cite{konstantinides:passalidis:2024d}.

We recall that a distribution $V$ is regularly varying with index $\alpha \in (0,\,\infty)$, denoted by $V \in \mathcal{R}_{-\alpha}$, if 
\beao
\lim \dfrac {\bV(t\,x)}{\bV(x)} = t^{-\alpha}\,,
\eeao
for any $t>0$. 

A distribution $F$ of a random vector ${\bf X}$ is (standard) multivariate regularly varying if there exists a univariate distribution $V \in \mathcal{R}_{-\alpha}$ with $\alpha \in (0,\,\infty)$and a non-degenerate Radon measure $\mu$ such that
\beam \label{eq.KLP.1.12}
\lim \dfrac 1{\bV(x)}\PP[{\bf X} \in x\,\bbb] = \mu(\bbb)\,,
\eeam
for every $\mu$-continuous Borel set $\bbb \subseteq \overline{\bbr}_+^d \setminus \{{\bf  0}\}$, that means $\mu(\partial \bbb)=0$. Symbolically, we write $F \in \mathrm{MRV}(\alpha,\,\mu)$. 

According to \cite[Prop. 4.14]{samorodnitsky:sun:2016} and \cite[Prop. 3.1]{konstantinides:passalidis:2024h}, the following inclusions hold:
\beam \label{eq.KLP.1.13}
\mathrm{MRV} \subsetneq \mathcal{A}_\mathscr{R} \subsetneq \mathcal{S}_\mathscr{R} \subsetneq \mathcal{L}_\mathscr{R} \,,
\eeam
where $\mathrm{MRV}$ denotes the union of all $\mathrm{MRV}(\alpha,\,\mu)$ classes for $\alpha \in (0,\,\infty)$ and any Radon measure $\mu$. These inclusions remain valid when replacing  $\mathscr{R}$ with any proper subset $A \subsetneq \mathscr{R}$.

The MRV structure defined by \eqref{eq.KLP.1.12} corresponds to the standard MRV,  which implies that the marginal distribution functions share the same regular variation index (under the non-degeneracy assumption $\mu({\bf y}\;:\;y_i > 0)>0$ for any $i=1,\,\ldots,\,d$). The non-standard MRV framework accommodates vectors $\mathbf{X}$ whose components exhibit heterogeneous tail behavior, meaning they may have different regular variation indices. We refer to Example \ref{exam.KLP.34.4} in Section \ref{sec.KLP.34}, which covers the case when the components of $\mathbf{X}$ have heterogeneous tail behaviors---when considering a specific set $A$ in \eqref{eq.KLP.34.1}, those with heavier tails may dominate. For detailed discussions on standard and non-standard MRV, see \cite{resnick:2007}, and for applications of non-standard MRV, refer to \cite{tang:yang:2019} and \cite{yang:cheng:zhang:2025}.

\section{Main results} \label{sec.KLP.3}

This section presents the principal results of our study, starting with the primary assumption required for the risk model.

\begin{assumption} \label{ass.KLP.3.1}
The sequences $\{{\bf X}^{(i)}\,,\;i \in \bbn\}$ and $\{N(t)\,,\;t\geq 0 \}$ are independent. 
\end{assumption} 
Assumption \ref{ass.KLP.3.1} represents a standard framework in both univariate and multivariate risk models. However, it excludes the time-dependent risk models that provide an important framework in actuarial practice. For models with specific dependence structures between claims and inter-arrival times, see, e.g, \cite{li:tang:wu:2010}, \cite{li:2016}, and \cite{cheng:konstantinides:wang:2022}, among others. Within the context of the standard framework, our first result establishes locally uniform asymptotic expressions for the probability that the discounted aggregate claims enter a `rare-set' $x\,A$, where 
$A \subsetneq \mathscr{R}$. This holds under the condition that the joint claim distribution $F$ belongs to the class of multivariate subexponential distributions on $A \subsetneq \mathscr{R}$

\bth \label{th.KLP.3.1}
Let $A \subsetneq \mathscr{R}$ be a fixed proper subset. For the discounted aggregate claims in \eqref{eq.KLP.1.2}, if Assumption \ref{ass.KLP.3.1} holds and $F \in \mathcal{S}_A$, then \eqref{eq.KLP.1.4} holds uniformly for any $t \in \Lambda_T$, where $T \in \Lambda$ is a fixed constant.
\ethe

The next theorem extends \eqref{eq.KLP.1.4} to a globally uniform estimate for any $t\in\Lambda$ by restricting the claim vector distribution to $\mathcal{A}_A$ and imposing a technical condition on inter-arrival times.

\bth \label{th.KLP.3.2}
Let $A \subsetneq \mathscr{R}$ be a fixed proper subset. Consider the discounted aggregate claims in \eqref{eq.KLP.1.2} under Assumption \ref{ass.KLP.3.1}. If $F \in \mathcal{A}_A$, and there exists $\vep>0$ such that $\PP[\theta_1 > \vep] = 1$ and the interest force $r$ is positive, then \eqref{eq.KLP.1.4} holds uniformly for all $t \in \Lambda$
\ethe

\bre \label{rem.KLP.3.1}
Two forms of particular relevance to insurance practice that the set $A$ can assume are
\beam \label{eq.KLP.3.14}
A_1 &=& \{{\bf x}\;:\;x_i > u_i\,,\;\exists \;i=1,\,\ldots,\,d \}\,,\\[2mm] \label{eq.KLP.3.15}
A_2 &=& \{{\bf x}\;:\;\sum_{i=1}^d l_i\,x_i > u \}\,,
\eeam
where $u_i > 0$, $u > 0$, and $l_1, \ldots, l_d > 0$ are constants satisfying $\sum_{i=1}^d l_i = 1$.
Hence, if the set $A$ takes the form \eqref{eq.KLP.3.14}, the asymptotic estimate \eqref{eq.KLP.1.4} characterizes the probability that any one of the $d$-lines of business exceeds its corresponding initial capital. Conversely, if $A$ takes the form \eqref{eq.KLP.3.15}, \eqref{eq.KLP.1.4} characterizes the asymptotic probability that the sum across all $d$-lines of business exceeds the insurer's total initial capital.
\ere

\bre \label{rem.KLP.3.2}
Recent literature has addressed the estimation of relation \eqref{eq.KLP.1.4} in related modeling contexts.  For instance, under more restrictive distribution classes and c\'{a}dl\'{a}g processes for the returns of financial risks,  \cite{konstantinides:passalidis:2024j} provides a more general estimation of \eqref{eq.KLP.1.4}, while \cite{passalidis:2025k} establishes locally uniform estimates for \eqref{eq.KLP.1.4} in a nonstandard risk model with multivariate strong subexponential claim vectors.
\ere

\bre \label{rem.KLP.3.3}
We first observe that in the univariate case ($d=1$), with $A=(1,\,\infty)$, and $r > 0$, Theorems \ref{th.KLP.3.1} and \ref{th.KLP.3.2} reduce to \cite[Th. 2.1, Th. 2.2]{hao:tang:2008}, respectively. Following the approach in the aforementioned work, we require the condition $\PP[\theta_1 > \vep]=1$, for some $\vep> 0.$ This is a technical assumption that does not impact actuarial applications, as $\vep>0$ can be chosen arbitrarily small.  More precisely, this condition---used in \cite[Lem. 4.3]{hao:tang:2008}---arises solely from technical considerations in their proof. Furthermore, the requirement $r > 0$ in Theorem \ref{th.KLP.3.2} is necessary only to ensure uniformity over the entire set $\Lambda.$ Specifically, to establish uniformity over the whole $\Lambda$, we must analyze the infinite sums in relation \eqref{eq.KLP.5.28}, which represent discounted aggregate claims over an infinite time horizon. The condition $r > 0$, together with the fact that ${\tau_{i}, i\in\mathbb{N}}$ are renewal epochs, ensures that the distribution of the infinite sum in \eqref{eq.KLP.5.28} is non-defective, given that $F$ is also non-defective (i.e., has no mass at $\infty$).

Finally, we note that if in Theorems \ref{th.KLP.3.1} and \ref{th.KLP.3.2}, we replace the  classes $\mathcal{S}_A$, $\mathcal{A}_A$ with $\mathcal{S}_\mathscr{R}$, $\mathcal{A}_\mathscr{R}$, respectively, then relation \eqref{eq.KLP.1.4} holds uniformly for all $t \in \Lambda_T$ (or $t \in \Lambda$, correspondingly) and for any $A \subsetneq \mathscr{R}$.
\ere

As in the univariate case, the asymptotic behavior of discounted aggregate claims directly determines the ruin probability in multivariate risk models. Although various definitions of ruin probability exist in this framework, we focus on the probability of the surplus process first entering a ruin set $L,$ which satisfies the following key assumption.

\begin{assumption} \label{ass.KLP.3.2}
We assume that the ruin set $L$ is open, decreasing (that means the $-L$ is increasing), its complement $L^c$, is convex, and it holds ${\bf 0} \in \partial L$. Additionally, for any $k>0$ we assume that it holds $k\,L=L$. 
\end{assumption} 
 
\bre \label{rem.KLP.3.4}
Assumption \ref{ass.KLP.3.2}, introduced in \cite{hult:lindskog:2006} and later employed in \cite{samorodnitsky:sun:2016}, provides valuable insights into ruin sets and their relationship with sets $A \subsetneq \mathscr{R}$. Particularly, in the univariate case ($d=1$), the ruin set $L = (-\infty, 0)$ satisfies both Assumption \ref{ass.KLP.3.2} and the relation $A = 1 - L$ where $A = (1, \infty) \subsetneq \mathscr{R}$. This relationship extends to the multivariate case through $A = \mathbf{l} - L \subsetneq \mathscr{R}$, where $\mathbf{l}$ denotes the deterministic weights used by the insurer to allocate initial capital across $d$-lines of business. Consequently, the following ruin sets directly correspond to the rare-event sets in \eqref{eq.KLP.3.14} and \eqref{eq.KLP.3.15} via this mapping:
\beao
&&L_1=\{{\bf x}\;:\; x_i<0\,,\;\exists\;i=1,\,\ldots,\,d \}\,,\\[2mm]\label{eq.KLP.3.16b}\notag
&&L_2=\{{\bf x}\;:\; \sum_{i=1}^d x_i<0\}\,.
\eeao
Recall that the ruin probability is defined as the probability that the surplus process enters a specified ruin set. Thus, ruin probability for set $L_1$ represents the probability that at least one of the $d$ business lines has a negative surplus, denoted by $\psi_{\text{or}}$. Ruin probability for set $L_2$ represents the probability that the aggregate surplus across all $d$ business lines is negative, denoted by $\psi_{\text{sum}}$. While the approach using ruin sets $L$ (and consequently classes $\mathcal{S}_\mathscr{R}$ and $\mathcal{A}_\mathscr{R}$) yields sufficiently general results for various dependence structures and claim distributions, it fails to provide asymptotic expressions for $\psi_{\text{and}}$ and $\psi_{\text{sim}}$ due to the non-convexity of $A^c$ with $A$ defined as follows:
\beao
A =\{{\bf x}\;:\;x_i> u_i\,,\;\forall\; i=1,\,\ldots,\,d\}\,,
\eeao
as it comes from \cite[Rem. 2.2]{konstantinides:passalidis:2024g}. For more discussions about $\psi_{or}$, $\psi_{sum}$, $\psi_{sim}$ and $\psi_{and}$ we refer to \cite{cheng:yu:2019}, \cite{lu:yuan:2022}.
\ere

Hence, for the ruin sets $L$ satisfying Assumption \ref{ass.KLP.3.2}, the ruin probability over a finite time horizon is defined as
\beam \label{eq.KLP.3.16.5} 
\psi_{{\bf l},\,L}(x,\,t) = \PP\left[{\bf U}(s) \in L\;:\;\exists \; s\in[0,\,t]\right]\,,
\eeam
and since we accepted $x\,L=L$ and $A= ({\bf l} - L) \subsetneq \mathscr{R}$ we obtain that it holds
\beam \label{eq.KLP.3.17} \notag
\psi_{{\bf l},\,L}(x,\,t) &=&\PP\left[{\bf D}_r(s) -\int_0^s e^{-rz}\,{\bf p}(z) dz \in x\, ({\bf l} - L)\;:\;\exists \; s\in[0,\,t]\right]\\[2mm]
&=&\PP\left[{\bf D}_r(s) -\int_0^s e^{-rz}\,{\bf p}(z) dz \in x\,A\;:\;\exists \; s\in[0,\,t]\right]\,,
\eeam
where
\beao
\int_0^t e^{-rz}\,{\bf p}(z) dz =\left( 
\begin{array}{c}
\int_{0}^{t} e^{-rz}\,p_{1}(z) dz \\ 
\vdots \\ 
\int_{0}^{t} e^{-rz}\,p_{d}(z) dz 
\end{array} 
\right)\,,
\eeao
for $t\geq 0$. Equation \eqref{eq.KLP.3.17} shows that the asymptotic behavior of ruin probability depends directly on the behavior of discounted aggregate claims, as established in Theorems \ref{th.KLP.3.1} and \ref{th.KLP.3.2}.

We now establish local and global uniform asymptotic estimates for the finite-time ruin probability in the risk model \eqref{eq.KLP.1.3}.
\bco \label{cor.KLP.3.1}
Consider a fixed set $A=({\bf l} - L) \subsetneq \mathscr{R}$ within the risk model defined by \eqref{eq.KLP.1.3}.
\begin{enumerate}
\item[(i)]
Under the conditions of Theorem \ref{th.KLP.3.1}, it holds that
\beam \label{eq.KLP.3.18}
&&\psi_{{\bf l},\,L}(x,\,t) \sim \int_0^t \PP\left[{\bf X} \in x\,e^{rs}\, A\right]\,\lambda(ds)\,,
\eeam
uniformly for all $t \in \Lambda_T$, and every fixed $T \in \Lambda$.
\item[(ii)]
Under the conditions of Theorem \ref{th.KLP.3.2}, the relation \eqref{eq.KLP.3.18} holds uniformly for any $t \in \Lambda$.
\end{enumerate}
\eco

Note that replacing the classes $\mathcal{S}_A$ and $\mathcal{A}_A$, with $\mathcal{S}_\mathscr{R}$ and  $\mathcal{A}_\mathscr{R}$, respectively, preserves the uniform validity of relation \eqref{eq.KLP.3.18} for all $t \in \Lambda_T$ (or $t \in \Lambda$) and any ruin set $L$ satisfying Assumption
\ref{ass.KLP.3.2} with $A=({\bf l} - L) \subsetneq \mathscr{R}$. 

\bre \label{rem.KLP.3.5}
Specifically, under the conditions of Theorem \ref{th.KLP.3.2}, and Corollary \ref{cor.KLP.3.1}(ii), respectively, and with the restriction that $F \in \mathrm{MRV}(\alpha,\,\mu)$, for some $\alpha \in (0,\,\infty)$ (which implies $F_A \in \mathcal{R}_{-\alpha}$, for any $A \subsetneq \mathscr{R}$) we can set $t=\infty$ in expressions \eqref{eq.KLP.1.4} and \eqref{eq.KLP.3.18}, and we obtain
\beam \label{eq.KLP.r.1} \notag
\PP\left[{\bf D}_r(\infty) \in x\,A\right]&\sim& \psi_{{\bf l},\,L}(x,\,\infty)\\[2mm]\notag
&\sim & \PP[{\bf X} \in x\,A]\,\int_0^{\infty}e^{-\alpha\,r\,s}\,\lambda(ds)\,\\[2mm]
&\sim &\mu(A)\bV(x)\,\int_0^{\infty}e^{-\alpha\,r\,s}\,\lambda(ds)\,,
\eeam
where the last step follows from \eqref{eq.KLP.1.12}, and we note that $\mu(A)\in (0,\infty)$ was established in the proof of \cite[Prop.4.14]{samorodnitsky:sun:2016}. Here, $\psi_{{\bf l},\,L}(x,\,\infty)$ denotes the infinite-time ruin probability, defined in \eqref{eq.KLP.3.16.5} with $s\in [0,\,\infty]$, which means $s\geq 0$. Furthermore, if $\{N(t)\,,\;t\geq 0\}$ is a homogeneous Poisson process with intensity $\lambda \in (0,\,\infty)$ then \eqref{eq.KLP.r.1} simplifies to
\beao
\PP[{\bf D}_r(\infty) \in x\,A] \sim \psi_{{\bf l},\,L}(x,\,\infty)\sim \dfrac {\lambda}{\alpha\,r}\,\mu(A)\bV(x)\,.
\eeao
\ere

\section{Examples of distribution from $\mathcal{S}_A$ and $\mathcal{A}_A$} \label{sec.KLP.34}

We now proceed with the construction of distribution examples belonging to the family $\mathcal{B}_A$, where $\mathcal{B} \in \{\mathcal{S},\,\mathcal{A}\}$ for some set $A$ of the form \eqref{eq.KLP.3.15} with $u=1$. These examples depart from the MRV framework as defined in \eqref{eq.KLP.1.13}, and instead rely solely on conditions imposed on the marginal distributions and the dependence structure among components, making them readily verifiable in practice. It should be noted, however, that subexponential marginal distributions alone are insufficient to guarantee membership in the class $\mathcal{S}_A$, as illustrated in \cite[Exam. 2.1]{passalidis:2025k}. For sufficient conditions ensuring that a distribution belongs to $\mathcal{S}_{\mathscr{R}}$, we refer to \cite[Sec. 4]{samorodnitsky:sun:2016}.

We now consider sets of the form
\beam \label{eq.KLP.34.1}
A = \left\{ {\bf x}\,:\, \sum_{i=1}^d l_i\,x_i > 1\right\}\,,
\eeam
where $l_1,\,\ldots,\,l_d > 0$ and $\sum_{i=1}^d l_i=1$. Although the following examples can be readily extended to sets of the form \eqref{eq.KLP.3.15}, we will work with \eqref{eq.KLP.34.1} for simplicity. Furthermore, we focus on the two-dimensional case $(d=2)$; the extension to arbitrary dimensions $d\in \bbn$ is straightforward for all examples presented.

In this section, we consider the random vector ${\bf X}=(X_1,\,X_2)$ that has non-negative components with marginal distributions $F_1$ and $F_2$. For the set $A$ given in \eqref{eq.KLP.34.1}, we can choose $I_A =\{(l_1,\,l_2)\}$. Then it holds that
\beam \label{eq.KLP.34.2}
\bF_A(x) =\PP[{\bf X}\in x\,A]=\PP \left[ \sum_{i=1}^2 l_i\,X_i > x\right]\,.
\eeam
Relations \eqref{eq.KLP.34.1} and \eqref{eq.KLP.34.2} demonstrate that the multivariate distribution of ${\bf X}$ is directly linked to the distribution of (deterministic) weighted sums. Consequently, the first example holds intrinsic value by establishing a closure property for the classes $\mathcal{S}$ and $\mathcal{A}$ under sums of dependent random variables with regression dependence, as introduced by \cite{lehmann:1966}.

\bexam \label{exam.KLP.34.1}
Consider the set $A$ defined in \eqref{eq.KLP.34.1}, with parameters $(l_1,\,l_2) \in [a,\,b]^2$, with $0 < a \leq b < 1$. Assume the random vector ${\bf X}=(X_1,\,X_2)$ has marginal distributions $F_1,\,F_2 \in \mathcal{B} \in \{\mathcal{S},\,\mathcal{A}\} $, with its components $X_1,\,X_2$ exhibiting regression dependence. That is, there exist constants $x_0>0$ and $C>0$ such that
\beam \label{eq.KLP.34.3}
\PP[ X_{i}>x_i\;|\;X_j=x_j]\leq C\,\PP \left[X_i > x_i\right]\,,
\eeam
for any $1\leq i \neq j \leq 2$ and all  $\min\{x_i,\, x_j\} \geq x_0$. Furthermore, assume $F_1*F_2 \in \mathcal{S}$. Then by \cite[Th. 2]{wang:2011} we obtain
\beam \label{eq.KLP.34.4}
\bF_A(x) =\PP[ l_1\,X_{1}+l_2\,X_{2}>x]\sim \PP [l_1\,X_1 > x_i ] +\PP[l_2\,X_2> x]=: \bF_1'(x) +\bF_2'(x) \,,
\eeam 
where $F_i'$ denotes the distribution of $l_i\,X_i$, for $i=1,\,2$. Since independence is a special case of \eqref{eq.KLP.34.3}, \cite[Th. 2]{wang:2011} also implies
\beam \label{eq.KLP.34.5}
\overline{F_1'*F_2'}(x) \sim \bF_1'(x) +\bF_2'(x) \,.
\eeam 
However, from the fact that $F_1',\,F_2' \in \mathcal{B} \in\{\mathcal{S},\,\mathcal{A}\} $ and relation \eqref{eq.KLP.34.5}, we apply \cite[Cor. 1.1]{leipus:siaulys:2020} for class $\mathcal{S}$ and \cite[Cor. 3.1]{konstantinides:passalidis:2024d} for class $\mathcal{A}$ to conclude that $F_1'*F_2' \in \mathcal{B} \in \{\mathcal{S},\,\mathcal{A}\}.$ Furthermore, from \eqref{eq.KLP.34.4} and \eqref{eq.KLP.34.5}, we obtain $\bF_A(x) \sim \overline{F_1'*F_2'}(x)$. Hence, by the closure property of class $\mathcal{B}$ under strong tail equivalence (see \cite[Th. 3]{teugels:1975} for $\mathcal{S}$ and \cite[Rem. 2.1]{konstantinides:passalidis:2024d} for $\mathcal{A}$), we conclude that $F_A \in \mathcal{B} $, which implies $F \in \mathcal{B}_A $.
\eexam

\bre \label{rem.KLP.34.1}
The closure property of the distribution class $\mathcal{S}$under convolution has been extensively studied, with several papers dedicated to this topic. For example, \cite{leslie:1989} provides a well-known counterexample demonstrating that $\mathcal{S}$ is not closed under convolution, while \cite{leipus:siaulys:2020} establishes necessary and sufficient conditions for closure. However, the closure of $\mathcal{S}$ under sums of dependent random variables remains relatively unexplored, despite existing work on the asymptotic tail behavior of sums of dependent variables (see, e.g., \cite{ko:tang:2008}, \cite{geluk:tang:2009}, \cite{jiang:gao:wang:2014}), among others. In our previous example, we observe that if $F_1,\,F_2 \in \mathcal{B} \in \{\mathcal{S},\,\mathcal{A}\} $ and $F_1*F_2 \in \mathcal{S},$ then the sum exhibits a certain insensitivity to the dependence structure specified in \eqref{eq.KLP.34.3}, yielding $F_{X_1+X_2} \in  \mathcal{B}.$  This result can be verified by applying the methodology of Example \ref{exam.KLP.34.1} with $l_1 = l_2 = 1/2$. 
\ere

For the next example, we introduce the class of dominatedly varying distributions, denoted by $\mathcal{D}$. A distribution $V$ belongs to $\mathcal{D}$ if, for any (or, equivalently, for some) $b\in (0,\,1)$, it satisfies
\beao
\limsup \dfrac{\bV(b\,x)}{\bV(x)}< \infty \,.
\eeao
From \cite{goldie:1978}, it is known that $\mathcal{D} \cap \mathcal{L}\equiv \mathcal{D} \cap \mathcal{S} $. A subclass of $\mathcal{D} \cap \mathcal{L}$ is the class of consistently varying distributions $\mathcal{C}$, which comprises all distributions satisfying
\beao
\lim_{v\downarrow 1}\liminf \dfrac{\bV(v\,x)}{\bV(x)}= 1 \,.
\eeao 
Recall that the following inclusion chain is well-established (see, e.g., \cite[Ch. 2]{leipus:siaulys:konstantinides:2023}):
\beao
\mathcal{R}:= \bigcup_{\alpha>0} \mathcal{R}_{-\alpha} \subsetneq \mathcal{C}\cap \mathcal{P_D} \subsetneq \mathcal{D}\cap \mathcal{A} \subsetneq \mathcal{A} \,.
\eeao

In the next example, we consider marginal distributions from subclasses of $\mathcal{S}$ and $\mathcal{A}$ that are narrower classes than those in Example \ref{exam.KLP.34.1}, while employing a broader dependence structure than the one specified in \eqref{eq.KLP.34.3}. Specifically, we consider tail asymptotic independence, as introduced in \cite{geluk:tang:2009}.

\bexam \label{exam.KLP.34.2}
Consider the set $A$ defined in \eqref{eq.KLP.34.1}. Assume the random vector ${\bf X}=(X_1,\,X_2)$, has marginal distributions $F_1,\,F_2 \in \mathcal{B}$ where $\mathcal{B}\in\{\mathcal{D}\cap \mathcal{S},\,\mathcal{D}\cap \mathcal{A}\} $ and its components $X_1,\,X_2$ satisfy the tail asymptotic independence condition:
\beam \label{eq.KLP.34.9}
\PP[X_i > x_i \:|\: X_j>x_j] \to 0 \,,
\eeam 
as $\min\{x_i,\,x_j\} \to \infty$, for any $1\leq i \neq j \leq 2$. By \cite[Th. 2.3]{li:2013}, we obtain
\beao
\bF_A(x)=\PP[l_1\,X_1 + l_2\,X_2 >x] \sim \PP[l_1\,X_1 >x] + \PP[l_2\,X_2 >x]\,.
\eeao
This implies  $F_A  \in \mathcal{B}$ with $\mathcal{B}\subset\mathcal{S}$ (or, respectively, $\mathcal{B}\subset\mathcal{A}$) and consequently $F \in \mathcal{S}_A$, (or, respectively, $F \in \mathcal{A}_A$).
\eexam

In the next example, we consider marginal distributions from the classes $\mathcal{D}\cap \mathcal{L}$,  $\mathcal{D}\cap \mathcal{A}$, or alternatively, from $\mathcal{C}$ and $\mathcal{C}\cap \mathcal{P_D}$, respectively. To accommodate these cases, we employ quasi-asymptotic independence---a dependence structure introduced in \cite{chen:yuen:2009} that generalizes the condition in \eqref{eq.KLP.34.9}.

\bexam \label{exam.KLP.34.3} 
Consider the set $A$ defined in \eqref{eq.KLP.34.1}. Assume the random vector ${\bf X}=(X_1,\,X_2)$, has marginal distributions $F_1,\,F_2 \in \mathcal{B},$ where $\mathcal{B}\in \{\mathcal{C},\,\mathcal{C}\cap \mathcal{P_D}\}$ and satisfies the quasi-asymptotic independence condition:
\beao
\lim \dfrac{\PP[X_1 > x\,,\: X_2>x]}{  \bF_1(x) +\bF_2(x)} =0 \,.
\eeao
By \cite[Th. 3.2]{chen:yuen:2009}, we obtain
\beao
\bF_A(x)=\PP[l_1\,X_1 + l_2\,X_2 >x] \sim \PP[l_1\,X_1 >x] + \PP[l_2\,X_2 >x]\,,
\eeao
which implies $F_A  \in \mathcal{B} $. Hence, $F_A \in \mathcal{S}$ (or respectively, $F_A \in \mathcal{A}$), and we conclude $F \in \mathcal{S}_A$ (or respectively, $F\in \mathcal{A}_A $).
\eexam
Examples \ref{exam.KLP.34.1}--\ref{exam.KLP.34.3} employ various forms of asymptotic independence between components to construct the classes $\mathcal{S}_A$ and $\mathcal{A}_A $. We now consider the following case with arbitrary dependence between components, where one marginal distribution belongs to $\mathcal{C}$ while the other belongs to $\mathcal{D}$, with the tail of the latter being asymptotically negligible relative to the former.

\bexam \label{exam.KLP.34.4} 
Consider the set $A$ defined in \eqref{eq.KLP.34.1}. Assume the random vector  ${\bf X}=(X_1,\,X_2)$ has marginal distributions $F_1 \in \mathcal{B},$ where $ \mathcal{B}\in \{\mathcal{C},\,\mathcal{C}\cap \mathcal{P_D}\} $ and $F_2 \in \mathcal{D}$, with $\bF_2(x) = o[\bF_1(x)]$. Then $F_1' \in \mathcal{B} \in \{\mathcal{C},\,\mathcal{C}\cap \mathcal{P_D}\} \subsetneq \mathcal{D}$ and $F_2' \in \mathcal{D}$ and moreover $\bF_2'(x) = o[\bF_1'(x)]$, due to $F_1,\,F_2\in\mathcal{D}$. Hence, by \cite[Lem. 3.3(i)]{yang:yuen:liu:2018} we obtain
\beao
\bF_A(x) = \PP[l_1\,X_1 + l_2\,X_2 > x] \sim \PP[l_1\,X_1 > x] = \bF_1'(x)\,.
\eeao
Therefore,  $F_A \in \mathcal{B} \in \{\mathcal{C},\,\mathcal{C}\cap \mathcal{P_D}\}$, we conclude that $F_A \in \mathcal{S}$, (or respectively, $F_A \in \mathcal{A}$), which implies $F \in \mathcal{S}_A$ (or respectively, $F \in \mathcal{A}_A$).
\eexam

\bre \label{rem.KLP.4.2}
Examples \ref{exam.KLP.34.1} - \ref{exam.KLP.34.4} yield concrete asymptotic expressions for relation \eqref{eq.KLP.1.4}. Under the assumptions of Examples \ref{exam.KLP.34.1} - \ref{exam.KLP.34.3}, Theorem \ref{th.KLP.3.1} combined with \eqref{eq.KLP.34.1} and \eqref{eq.KLP.34.4} provides the uniform asymptotic estimate
\beam \label{eq.KLP.34.a} \notag
\PP\left[{\bf D}_r(t) \in x\,A \right]&=&\PP\left[l_1 \sum_{i=1}^{N(t)} X_1^{(i)}\,e^{-r\,\tau_i} + l_2 \sum_{i=1}^{N(t)} X_2^{(i)}\,e^{-r\,\tau_i} >x\right]\\[2mm] 
&\sim & \int_0^t \left( \PP\left[l_1 \,X_1>x\,e^{r\,s}\right] + \PP\left[l_2 \,X_2>x\,e^{r\,s}\right] \right)\,\lambda(ds)\,,
\eeam
for all $t \in \Lambda_T$. Furthermore, if $r>0$, $\PP[\theta_1 > \vep]=1$ for some $\vep>0$, and distributions from class $\mathcal{A}$ are used instead of $\mathcal{S}$ in Examples \ref{exam.KLP.34.1}–\ref{exam.KLP.34.3}, then \eqref{eq.KLP.34.a} holds uniformly over $t \in \Lambda$.

Next, we adopt the assumptions of Example \ref{exam.KLP.34.4}. By Theorem \ref{th.KLP.3.1}, we have $\PP[Y_A > x] \sim \PP[l_1,\,X_1 > x]$ (or in the case $r>0$, $\PP[\theta_1>\vep]=1$ and marginals from $\mathcal{C}\cap \mathcal{P_D}$) which yields the asymptotic expression:
\beam \label{eq.KLP.34.b} \notag
\PP\left[{\bf D}_r(t) \in x\,A \right]&=&\PP\left[l_1 \sum_{i=1}^{N(t)} X_1^{(i)}\,e^{-r\,\tau_i} + l_2 \sum_{i=1}^{N(t)} X_2^{(i)}\,e^{-r\,\tau_i} >x\right]\\[2mm]
&\sim&\int_0^t \PP\left[l_1 \,X_1>x\,e^{r\,s}\right]\,\lambda(ds)\,,
\eeam
uniformly for all $t \in \Lambda_T$, or $t \in \Lambda$, respectively.

We note that by Corollary \ref{cor.KLP.3.1}, \eqref{eq.KLP.34.a} and \eqref{eq.KLP.34.b} also hold for the ruin probabilities $\psi_{{\bf l},\,L}$, where $A=({\bf l} - L)$ is given by \eqref{eq.KLP.34.1} with $d=2.$ Furthermore, this ruin probability is closely related to $\psi_{sum},$ though more general due to the coefficients $l_1,\,l_2.$ The expressions in \eqref{eq.KLP.34.a} and \eqref{eq.KLP.34.b} represent standard formulas for bivariate risk models; see, for example, \cite{chen:wang:wang:2013}, \cite{gao:yang:2014}, and \cite{chen:li:cheng:2023}.
\ere

In Examples \ref{exam.KLP.34.1}–\ref{exam.KLP.34.3}, we employed a specific form of asymptotic independence between components. Example \ref{exam.KLP.34.4} considered arbitrary dependence, with one tail dominating the other. We now examine a case of asymptotic dependence with strongly equivalent tails. However, instead of the set in \eqref{eq.KLP.34.1}, we restrict to
\beam \label{eq.KLP.E.1}
A = \left\{ {\bf x}\;:\; \sum_{i=1}^d \dfrac 1d x_i >1 \right\}\,,
\eeam
and for simplicity, consider only $d=2$. Thus, we study 
\beao
\bF_A(x) = \PP\left[\dfrac{X_1}2 + \dfrac{X_2}2 > x\right]\,,
\eeao  
for the set defined in \eqref{eq.KLP.E.1}. We utilize a multivariate Gumbel distribution, extensively studied in \cite{kluppelberg:resnick:2008}, \cite{asimit:furman:tang:vernic:2011}, and \cite{li:2022a}. Within our framework, we note that if a distribution $V$ belongs to the Gumbel maximum domain of attraction, denoted by $V \in \mathrm{MDA}(\Lambda)$ (see, e.g., \cite{haan:ferreira:2006}), then $V \in \mathcal{R}_{-\infty}.$ This implies \beao
\lim \dfrac{\bV(v\,x)}{\bV(x)}=0\,,
\eeao
for any $v>1$. Consequently, $\mathcal{R}_{-\infty} \subsetneq \mathcal{P_D}$, and therefore $\mathcal{S}\cap \mathcal{R}_{-\infty} \subsetneq \mathcal{A}$.

\bexam \label{exam.KLP.4.5} 
Consider the set $A$ defined in \eqref{eq.KLP.E.1}. Following Assumption 2.2 of \cite{asimit:furman:tang:vernic:2011}, assume the distribution $F$ of the non-negative random vector ${\bf X}=(X_1,\,X_2)$ satisfies
\beam \label{eq.KLP.E.2}
\lim \dfrac{\PP[X_1 > x+a(x)\,z_1\,,\;X_2 > x+a(x)\,z_2]}{\bF(x)}=H_G({\bf z})\,,
\eeam
for all ${\bf z} \in \bbr^2,$ where $H_G$ is non-degenerate and $a(\cdot)\;:\;[0,\,\infty) \to (0,\,\infty)$ is an auxiliary function satisfying  $a(x)=o(x)$ and $a(x+a(x)\,z) \sim a(x)$ locally uniformly in $z$ (see \cite[Ch. 1.2]{haan:ferreira:2006}). We further assume $F_1 \in \mathrm{MDA}(\Lambda) \cap \mathcal{S}$. From \eqref{eq.KLP.E.2} and the argument in \cite[Sec. 2.2]{asimit:furman:tang:vernic:2011}, we obtain
\beao
\lim \dfrac{\PP\left[\left(\dfrac{X_1 -x}{a(x)}\,,\;\dfrac{X_2 -x}{a(x)} \right) \in \bbb \right]}{\bF_1(x)}=\mu(\bbb)\,,
\eeao
for any Borel set $\bbb \in [-\infty\,,\;\infty]^2 \setminus \{-(\infty\,,\;\infty)\}$, where
\beam \label{eq.KLP.E.3}
\mu [(z_1,\;\infty)\times (z_2,\;\infty)]=H_G({\bf z})>0\,.
\eeam
This implies that
\beao
\lim \dfrac{\bF_2(x)}{\bF_1(x)}=\mu([-\infty,\,\infty]\times(0,\,\infty])> 0\,,
\eeao
establishing tail equivalence between $F_1$ and $F_2.$ Consequently, it holds $F_1,\,F_2 \in \mathcal{S} \cap \mathcal{R}_{-\infty} \subsetneq \mathcal{A} $. From \eqref{eq.KLP.E.2} and \cite[Cor. 4.1]{kluppelberg:resnick:2008}, we have
\beam \label{eq.KLP.E.4}
\PP[X_1 + X_2 >2x] \sim \mu \left({\bf x}\;:\;\sum_{i=1}^2 x_i >0\right)\,\bF_1(x)\,
\eeam
see \cite[Cor 4.1]{kluppelberg:resnick:2008}. Thus, by \eqref{eq.KLP.E.4},
\beao
\bF_A(x)=\PP\left[\dfrac{X_1}2 + \dfrac{X_2}2 >x \right] \sim \mu \left({\bf x}\;:\;\sum_{i=1}^2 x_i >0 \right)\,\bF_1(x)\,,
\eeao
since $F_1 \in \mathcal{S}\cap \mathcal{R}_{-\infty} \subsetneq \mathcal{A} $ and \eqref{eq.KLP.E.3} imply
\beao
 \mu \left({\bf x}\;:\;\sum_{i=1}^2 x_i >0 \right) >0\,,
\eeao
the closure property of $\mathcal{A}$ with tail equivalence implies $F_A \in \mathcal{A}.$ Hence, $F \in \mathcal{A}_A$.
\eexam

We conclude this section with an example of a distribution in $\mathcal{A}_{\mathscr{R}}$, inspired by the construction method for $\mathcal{S}_{\mathscr{R}}$ using rotationally invariant random variables, as presented in \cite[Exam. 4.17]{samorodnitsky:sun:2016}. Let ${\bf X}=(X_1,\,\ldots,\,X_d)$, following distribution $F$, be a non-negative random vector with distribution $F$ such that $\PP[{\bf X}={\bf 0}]=0$, and let $\Delta_d \subset \bbr^d$ denote the unit sphere.

\bexam \label{exam.KLP.34.5}
Let $A \subsetneq \mathscr{R}$ be a fixed set. Following the construction in \cite[Exam. 4.17]{samorodnitsky:sun:2016}, we define a one-dimensional random variable $Z$ with distribution $V$, whose tail distribution is given by
\beao
\bV(x)=\PP\left[\sum_{i=1}^d X_i > x\:\Big|\:\dfrac{{\bf X}}{\sum_{i=1}^d X_i} \in (\theta_1,\,\ldots,\,\theta_d)\right]\,,
\eeao 
for any $ (\theta_1,\,\ldots,\,\theta_d) \in \Delta_d$. We also define the random variable
\beam \label{eq.KLP.34.12}
H=\inf \left\{ u>0\::\: u\,\dfrac{{\bf X}}{\sum_{i=1}^d X_i} \in A\right\}\,,
\eeam 
for any $A\subsetneq \mathscr{R}$. If $Z$ and $H$ are independent, we define a rotationally invariant distribution $F_A$, for the random variable $Y_A$, where $Y_A$ is given by
\beam \label{eq.KLP.34.13}
Y_A\stackrel{d}{=} Z\,H^{-1}\,,
\eeam
where `$\stackrel{d}{=} $' denotes equality in distribution.

 By \eqref{eq.KLP.34.12}, $H$ is bounded away from zero; hence, $H^{-1}$ is bounded from above. Therefore, if $V\in\mathcal{A}$, then by \eqref{eq.KLP.34.13} and the independence of $Z$ and $H$, we can apply 
\cite[Theorem 2.1]{tang:2006} to conclude that $F_A  \in \mathcal{A}$. Hence, for any $A\subsetneq \mathscr{R},$ we obtain $F \in \mathcal{A}_A$, which implies $F \in \mathcal{A}_{\mathscr{R}}$.
\eexam

\section{Proof of Theorem \ref{th.KLP.3.1}} \label{sec.KLP.4}

We now prove Theorem \ref{th.KLP.3.1} using a preliminary lemma. The following lemma generalizes \cite[Lem. 6.1]{passalidis:2025k} under conditions that are more general than those required in Theorem \ref{th.KLP.3.1}. In what follows, we use the notation
\beao
&&F_i(x\,A)=\PP\left[ {\bf X}^{(i)} \in x\,A\right]\,,\\[2mm] 
&&\bF_A^{(i)}(x)=\PP\left[Y^{(i)}_A>x \right]\,,
\eeao 
with 
\beao
Y^{(i)}_A=\sup\left\{u\;:\;  {\bf X}^{(i)} \in u\,A\right\}\,,
\eeao 
for any $i=1,\,\ldots,\,n$. 

\ble \label{lem.KLP.4.1}
Let $A\subsetneq \mathscr{R}$ be a fixed set and consider $n$ independent, non-negative random vectors $ {\bf X}^{(1)},\,\ldots,\,{\bf X}^{(n)} $ with distributions $F_1,\,\ldots,\,F_n \in \mathcal{L}_A$. We assume there exists a distribution $F \in \mathcal{S}_A$ such that 
\beao
F_i(x\,A)\asymp F(x\,A)\,,
\eeao 
for any $i=1,\,\ldots,\,n$. Then, for any fixed $0 <a \leq b < \infty$, it holds
\beam \label{eq.KLP.4.19}
\PP\left[\sum_{i=1}^n c_i\,{\bf X}^{(i)} \in x\,A \right] \sim \sum_{i=1}^n \PP\left[c_i\,{\bf X}^{(i)} \in x\,A\right]\,,
\eeam
uniformly for any ${\bf c}_n:=(c_1,\,\ldots,\,c_n) \in [a,\,b]^n$.
\ele

\pr~
We begin by establishing the upper bound for \eqref{eq.KLP.4.19}. Applying \cite[Prop. 2.4]{konstantinides:passalidis:2024g} followed by \cite[Lem. 1]{tang:yuan:2014}, and noting that $F_A^{(1)},\,\ldots,\,F_A^{(n)} \in \mathcal{L}$, $F_A \in \mathcal{S}$ and 
\beao
\bF_A^{(i)}(x)\asymp \bF_A(x)\,,
\eeao 
for all $i=1,\,\ldots,\,n$, we obtain
\beam \label{eq.KLP.4.20} \notag
\PP\left[\sum_{i=1}^n c_i {\bf X}^{(i)} \in x A \right] &\leq &\PP\left[\sum_{i=1}^n c_i Y_A^{(i)} > x \right] \\[2mm]\notag
&\sim &\sum_{i=1}^n \PP\left[c_i Y_A^{(i)} > x\right]\\[2mm]&=&\sum_{i=1}^n \PP\left[c_i{\bf X}^{(i)} \in x A \right] ,
\eeam
uniformly, for any ${\bf c}_n \in [a,\,b]^n$.

For the lower bound, since $A$ is an increasing set, and the random vectors ${\bf X}^{(i)}$ are non-negative and independent, the Bonferroni inequality yields
\beam \label{eq.KLP.4.21} \notag
\PP\left[\sum_{i=1}^n c_i\,{\bf X}^{(i)} \in x\,A \right] &\geq&\PP\left[\bigcup_{i=1}^n \{c_i\,\,{\bf X}^{(i)} \in x\,A\}\right] \\[2mm] \notag
&\geq& \sum_{i=1}^n\PP\left[ c_i\,{\bf X}^{(i)} \in x\,A \right] - \sum_{1\leq i<j \leq n} \PP\left[c_i\,{\bf X}^{(i)} \in x\,A\,,\;c_j\,{\bf X}^{(j)} \in x\,A \right]\\[2mm] \notag
&=&\sum_{i=1}^n \PP\left[c_i\,Y_A^{(i)} > x\right]- \sum_{1\leq i<j \leq n} \PP\left[c_i\,Y_A^{(i)} > x\,,\;c_j\,Y_A^{(j)} > x \right]\\[2mm]\notag
&\sim &\sum_{i=1}^n \PP\left[c_i\,Y_A^{(i)} > x\right]\\[2mm]&=& \sum_{i=1}^n \PP\left[ c_i\,{\bf X}^{(i)} \in x\,A \right]\,,
\eeam
uniformly for any ${\bf c}_n \in [a,\,b]^n$. Hence, by relations \eqref{eq.KLP.4.20} and \eqref{eq.KLP.4.21} we get the \eqref{eq.KLP.4.19}.

This completes the proofs of Theorem \ref{th.KLP.3.1}.
~\halmos

{\bf Proof of Theorem \ref{th.KLP.3.1}}~
Let $M \in \bbn$ be an integer and fix $t \in \Lambda_T$. By the total probability theorem,
\beam \label{eq.KLP.4.22} \notag
\PP\left[{\bf D}_r(t) \in x\,A \right] &=& \PP\left[\sum_{i=1}^{N(t)}  {\bf X}^{(i)}\,e^{-r \tau_i} \in x\,A \right] \\[2mm] \notag
&=&\left( \sum_{n=1}^M + \sum_{n=M+1}^{\infty} \right)\PP\left[ \sum_{i=1}^{n}{\bf X}^{(i)}\,e^{-r \tau_i} \in x\,A\,,\;N(t)=n \right]\\[2mm] 
&=:&I_1(x,\,t,\,M)+I_2(x,\,t,\,M)\,,
\eeam
We start with the estimation of $I_2(x,\,t,\,M)$. For a fixed $\delta > 0$, the following holds:
\beao
I_2(x,\,t,\,M) &\leq &\sum_{n=M+1}^{\infty}\PP\left[ \sum_{i=1}^{n}{\bf X}^{(i)}\,e^{-r \tau_1} \in x\,A\,,\;\tau_n \leq t \right]\\[2mm]
&&= \sum_{n=M+1}^{\infty} \int_0^t \PP\left[ \sum_{i=1}^{n}{\bf X}^{(i)}\,e^{-r\,s} \in x\,A \right]\,\PP[\tau_{n-1} \leq t-s]\,\PP[\tau_1 \in ds]\\[2mm]
&&=\sum_{n=M+1}^{\infty} \int_0^t \PP\left[ \sum_{i=1}^{n}{\bf X}^{(i)}\,e^{-r\,s} \in x\,A \right]\,\PP[N(t-s) \geq n-1]\,\PP[\tau_1 \in ds]\\[2mm]
&&\leq k\sum_{n=M+1}^{\infty}\,\PP[N(T) \geq n-1](1+\delta)^n \int_0^t \PP\left[ {\bf X}\,e^{-r s} \in x\,A \right]\,\PP[\tau_1 \in ds]\\[2mm]
&&\leq k\left(\int_0^t \PP\left[ {\bf X}\,e^{-r s} \in x\,A \right]\,\lambda(ds) \right) \sum_{n=M+1}^{\infty}\,\PP[N(T) \geq n-1](1+\delta)^n \,,
\eeao
uniformly for all $t \in \Lambda_T$, where the constant $k>0$ exists for any  $\delta > 0$ by \cite[Prop. 4.12(c)]{samorodnitsky:sun:2016} (recall that $F\in \mathcal{S}_A$). Since $\{N(t)\,,\;t \geq 0\}$ is a renewal process and $T<\infty$, it has an analytic renewal function near zero (see \cite{stein:1946}). Therefore, for any  $\delta' \in (0,\,1),$ there exists a sufficiently large $M$ such that 
\beam \label{eq.KLP.4.23} 
&&I_2(x,\,t,\,M) \leq k\,\delta'\,\int_0^t \PP\left[ {\bf X} \in x\,e^{r s}\,A \right]\,\lambda(ds) \,,
\eeam
uniformly, for all $t \in \Lambda_T$.

For $I_1(x,\,t,\,M)$, we have
\beam \label{eq.KLP.4.24} \notag
&&I_1(x,\,t,\,M) \\[2mm] \notag
&=&\sum_{n=1}^M \iint_{\{0\leq \tau_1 \leq \cdots \leq  \tau_n \leq t <\tau_{n+1}\}} \PP\left[ \sum_{i=1}^n {\bf X}^{(i)}\,e^{-r s_i} \in x\,A \right]\,\PP[\tau_1 \in ds_1,\,\ldots,\,\tau_{n+1} \in ds_{n+1}] \\[2mm] \notag
&\sim&\sum_{n=1}^M \iint_{\{0\leq \tau_1 \leq \cdots \leq \tau_n \leq t <\tau_{n+1}\}} \sum_{i=1}^n \PP\left[ {\bf X}^{(i)}\,e^{-r s_i} \in x\,A \right]\,\PP[\tau_1 \in ds_1,\,\ldots,\,\tau_{n+1} \in ds_{n+1}]\\[2mm]\notag
&=&\sum_{n=1}^M \,\sum_{i=1}^n \PP\left[ {\bf X}^{(i)}\,e^{-r \tau_i} \in x\,A\,,\;N(t)=n \right]\\[2mm]\notag
&=&\left(\sum_{n=1}^{\infty} - \sum_{n=M+1}^{\infty} \right)\sum_{i=1}^n \PP\left[ {\bf X}^{(i)}\,e^{-r \tau_i} \in x\,A\,,\;N(t)=n \right]\\[2mm]
&=:&I_{11}(x,\,t,\,M) - I_{12}(x,\,t,\,M) \,,
\eeam
uniformly for all $t \in \Lambda_T$. The asymptotic equivalence in the second step follows from Lemma \ref{lem.KLP.4.1} via the dominated convergence theorem, applicable by the generalized multivariate Kesten's Lemma in \cite[Lem. 5.2]{konstantinides:passalidis:2024g}.

For $I_{11}(x,\,t,\,M)$, interchanging summations yields
\beam \label{eq.KLP.4.25} 
I_{11}(x,\,t,\,M) \sim \sum_{i=1}^{\infty} \PP\left[ {\bf X}^{(i)}\,e^{-r \tau_i} \in x\,A\,,\;\tau_i \leq t \right]=\int_0^t \PP\left[ {\bf X} \in x\,e^{r s}\,A \right]\,\lambda(ds) \,,
\eeam
uniformly for all $t \in \Lambda_T$. For $I_{12}(x,\,t,\,M)$, 
\beao
I_{12}(x,\,t,\,M) &\leq& \sum_{n=M+1}^{\infty} \sum_{i=1}^{n} \PP\left[ {\bf X}^{(i)}\,e^{-r \tau_1} \in x\,A\,,\;\tau_n \leq t \right]\\[2mm]
&=&\sum_{n=M+1}^{\infty} \sum_{i=1}^{n} \int_0^t \PP\left[ {\bf X}\,e^{-r s} \in x\,A \right]\,\PP[\tau_{n-1} \leq t-s ]\,\PP[\tau_1 \in ds]\\[2mm]
&\leq& \left(\sum_{n=M+1}^{\infty} n \PP[N(T) \geq n-1 ]\right)\,\int_0^t \PP\left[ {\bf X} \in x\,e^{r s}\,A \right]\,\lambda(ds) \,,
\eeao
holds uniformly for all $t \in \Lambda_T$. Since $\{N(t)\,,\;t\geq 0\}$ is a renewal counting process and $T<\infty,$ for any $\delta' \in (0,\,1),$ there exists a sufficiently large $M,$ such that 
\beam \label{eq.KLP.4.26} 
I_{12}(x,\,t,\,M) \leq \delta' \int_0^t \PP\left[ {\bf X} \in x\,e^{r s}\,A \right]\,\lambda(ds) \,,
\eeam
holds uniformly for all $t \in \Lambda_T$. From \eqref{eq.KLP.4.24}, \eqref{eq.KLP.4.25}, and \eqref{eq.KLP.4.26}, we obtain
\beam \label{eq.KLP.4.27}
(1- \delta')\int_0^t \PP\left[ {\bf X} \in x\,e^{r s}\,A \right]\,\lambda(ds)\lesssim I_1(x,\,t,\,M) \lesssim \int_0^t \PP\left[ {\bf X} \in x\,e^{r s}\,A \right]\,\lambda(ds) \,,
\eeam
holds uniformly for any $t \in \Lambda_T$. Hence, substituting relations \eqref{eq.KLP.4.23}, \eqref{eq.KLP.4.27} into \eqref{eq.KLP.4.22}, and given the arbitrariness of $\delta'>0$, we conclude that relation \eqref{eq.KLP.1.4} holds uniformly for any $t \in \Lambda_T$.
~\halmos

\section{Proofs of Theorem \ref{th.KLP.3.2} and Corollary \ref{cor.KLP.3.1} } \label{sec.KLP.5}

We now proceed to the proof of Theorem \ref{th.KLP.3.2} and Corollary \ref{cor.KLP.3.1}, beginning with a preliminary lemma. This result extends \cite[Lem. 4.3]{hao:tang:2008} to the multivariate setting. Define
\beam \label{eq.KLP.5.28}
{\bf D}_r(\infty):=\sum_{i=1}^{\infty} {\bf X}^{(i)}\,e^{-r \tau_i}\,.
\eeam

\ble \label{lem.KLP.5.1}
Let $A\subsetneq\mathscr{R}$ be a fixed set. Under the assumptions of Theorem \ref{th.KLP.3.2}, the following inequality holds:
\beao
\PP\left[ {\bf D}_r(\infty) \in x\,A \right] \lesssim \int_0^{\infty} \PP\left[ {\bf X}  \in x\,e^{r s}\,A \right]\,\lambda(ds)\,.
\eeao
\ele

\pr~
Since $Y_A^{(i)}$ follows distribution $F_A \in \mathcal{A}$, are mutually independent, $r>0$ and $\PP[\theta_1>\vep]=1$, for some $\vep>0$, we obtain that it holds
\beao
\PP\left[ {\bf D}_r(\infty) \in x\,A \right] &=& \PP\left[ \sum_{i=1}^{\infty} {\bf X}^{(i)}\,e^{-r \tau_i} \in x\,A \right]\\[2mm]&=&\PP\left[\sup_{{\bf p} \in I_A} {\bf p}^T\left(\sum_{i=1}^{\infty} {\bf X}^{(i)}\,e^{-r \tau_i} \right)> x \right] \\[2mm]\notag
&\leq& \PP\left[\sum_{i=1}^{\infty} Y_A^{(i)}\,e^{-r \tau_i} > x \right]  \\[2mm]\notag&\lesssim& \int_0^{\infty} \bF_A(x\,e^{r s})\,\lambda(ds)\\[2mm]&=&\int_0^{\infty} \PP\left[ {\bf X}  \in x\,e^{r s}\,A \right]\,\lambda(ds)\,,
\eeao
where the fourth step employs \cite[Lem. 4.3]{hao:tang:2008}, concluding the result.

~\halmos

Now, we proceed with the proof of Theorem \ref{th.KLP.3.2}.

{\bf Proof of Theorem \ref{th.KLP.3.2}}~
By \cite[Lem. 4.2]{hao:tang:2008}, for any $\delta>0$, there exist sufficiently large $T' \in \Lambda$ and $x_0>0$, such that 
\beam \label{eq.KLP.5.30} \notag
\int_{T'}^{\infty} \PP\left[ {\bf X}  \in x\,e^{r s}\,A \right]\,\lambda(ds)&=&\int_{T'}^{\infty} \PP\left[ Y_A > x\,e^{r s} \right]\,\lambda(ds)\\[2mm]\notag
& \leq &\delta\,\int_0^{T'} \PP\left[ Y_A > x\,e^{r s} \right]\,\lambda(ds)\\[2mm]&=&\delta\,\int_0^{T'} \PP\left[  {\bf X}  \in x\,e^{r s}\,A \right]\,\lambda(ds)\,,
\eeam
for all $x\geq x_0$. For $t \in (T',\,\infty]$, Theorem \ref{th.KLP.3.1} and \eqref{eq.KLP.5.30} yield
\beam \label{eq.KLP.5.31} \notag
\PP\left[ {\bf D}_r(t) \in x\,A \right]
&\geq &\PP\left[ {\bf D}_r(T') \in x\,A \right] \\[2mm]\notag
&\sim &\int_0^{T'} \PP\left[ {\bf X}  \in x\,e^{r s}\,A \right]\,\lambda(ds)\\[2mm]\notag
&\geq & \left(\int_0^{t} -\int_{T'}^{\infty}\right)  \PP\left[ {\bf X}  \in x\,e^{r s}\,A \right]\,\lambda(ds) \\[2mm]
&\geq& (1-\delta)\,\int_0^{t} \PP\left[  {\bf X}  \in x\,e^{r s}\,A \right]\,\lambda(ds)\,,
\eeam
uniformly for all $t \in (T',\,\infty]$.

Conversely, Lemma \ref{lem.KLP.5.1} and \eqref{eq.KLP.5.30} imply
\beam \label{eq.KLP.5.32} \notag
\PP\left[ {\bf D}_r(t) \in x\,A \right] &\leq& \PP\left[ {\bf D}_r(\infty) \in x\,A \right]\\[2mm]   \notag
&\lesssim& \int_0^{\infty} \PP\left[ {\bf X}  \in x\,e^{r s}\,A \right]\,\lambda(ds)\\[2mm]\notag
&\leq &\left(\int_0^{t} +\int_{T'}^{\infty}\right)  \PP\left[ {\bf X}  \in x\,e^{r s}\,A \right]\,\lambda(ds) \\[2mm]&\leq &(1+\delta)\,\int_0^{t} \PP\left[  {\bf X}  \in x\,e^{r s}\,A \right]\,\lambda(ds)\,,
\eeam
uniformly for all  $t \in (T',\,\infty]$. By \eqref{eq.KLP.5.31} and \eqref{eq.KLP.5.32}, and the arbitrariness of $\delta>0$, \eqref{eq.KLP.1.4} holds uniformly for all $t \in (T',\,\infty]$. Theorem \ref{th.KLP.3.1} establishes uniform convergence for all $t \in \Lambda_{T'},$ so \eqref{eq.KLP.1.4} holds uniformly for all $t \in \Lambda$.   
~\halmos

{\bf Proof of Corollary \ref{cor.KLP.3.1}}~
We establish part (ii), as part (i) follows analogously by substituting Theorem \ref{th.KLP.3.1} for Theorem \ref{th.KLP.3.2}.
The premiums have bounded densities, with $C_i\in (0,\,\infty)$ bounding the $i$-th premium density. Thus,
\beam \label{eq.KLP.5.33}
0 \leq \int_0^{t} \,e^{-r y}\,p_i(y) \,dy \leq C_i\,t < \infty\,,
\eeam
for all $t \in \Lambda$, and $i=1,\,\ldots,\,d.$ Given $A=({\bf l} -L) \subsetneq \mathscr{R}$ and \eqref{eq.KLP.5.33}, \cite[Lem. 4.3(d)]{samorodnitsky:sun:2016} implies that for any ${\bf a} \in \bbr^d$, there exists $K>0$, such that for all $x> K,$ 
\beao
(x+K)\,A \subsetneq x\,A + {\bf a} \subsetneq (x-K)\,A\,.
\eeao   
Consequently,
\beam \label{eq.KLP.5.34} \notag
\psi_{{\bf l},L}(x,\,t) &\leq& \PP\left[ {\bf D}_r(t) -\int_0^{s} \,e^{-r y}\,{\bf p} (y) \,dy\in x\,A\,,\; \exists \;s\in[0,\,t] \right] \\[2mm] \notag
&\leq&  \PP\left[ {\bf D}_r(t) \in (x-K)\,A\right]\\[2mm] \notag &\sim &\int_0^t \PP\left[ {\bf X}  \in (x-K)\,e^{r s}\,A \right]\,\lambda(ds) \\[2mm] 
&\sim&\int_0^t \PP\left[ {\bf X}  \in x\,e^{r s}\,A \right]\,\lambda(ds)  \,,
\eeam
uniformly for all $t \in \Lambda.$ The third step applies Theorem \ref{th.KLP.3.2}, and the last step uses the uniformity property of class  $\mathcal{L}$ (since $F_A \in \mathcal{A} \subsetneq \mathcal{L}$), as established in \cite[Lem. 3.1]{cheng:cheng:2018}.

Similarly,
\beam \label{eq.KLP.5.35}\notag
\psi_{{\bf l},L}(x,\,t) &\geq& \PP\left[ {\bf D}_r(t) -\int_0^{t} \,e^{-r y}\,{\bf p}(y) \,dy\in x\,A \right] \\[2mm] \notag
&\gtrsim & \int_0^t \PP\left[ {\bf X}  \in (x+K)\,e^{r s}\,A \right]\,\lambda(ds)  \\[2mm]
&\sim &\int_0^t \PP\left[ {\bf X} \in x\,e^{r s}\,A \right]\,\lambda(ds)  \,,
\eeam
uniformly for all $t \in \Lambda$. From \eqref{eq.KLP.5.34} and \eqref{eq.KLP.5.35}, \eqref{eq.KLP.3.18} holds uniformly for all $t \in \Lambda$.
\halmos

\noindent \textbf{Acknowledgments.} 
Most of the paper was written during the visit of D.K. and C.P. to the Department of Financial and Actuarial Mathematics at Xi'an Jiaotong-Liverpool University. The research of J.L. was supported by the National Natural Science Foundation of
China (NSFC: 12201507; 72171055) and the XJTLU Research Enhancement Fund (REF-22-02-003).

We would like to thank the two anonymous referees for their comments, which improved the paper.

\noindent \textbf{Disclosure statement.}
There are no financial or non-financial competing interests.

\end{document}